\documentclass{svproc}
\usepackage{url}

\usepackage{amsmath}
\usepackage{graphicx}

\begin{document}
\mainmatter              % start of a contribution
\title{Parallel efficiency of monolithic and fixed-strain solution strategies for poroelasticity problems}
\titlerunning{Parallel efficiency for poroelasticity}  % abbreviated title (for running head)
\author{Denis Anuprienko} %\and Roger Temam\inst{2}
%Jeffrey Dean \and David Grove \and Craig Chambers \and Kim~B.~Bruce \and
%Elsa Bertino}
%
\authorrunning{D. Anuprienko} % abbreviated author list (for running head)
\tocauthor{Denis Anuprienko}
\institute{Nuclear Safety Institute RAS, Moscow 115191, Russian Federation
\email{anuprienko@ibrae.ac.ru}
}
\maketitle              % typeset the title of the contribution

\begin{abstract}
	Poroelasticity is an example of coupled processes which are crucial for many applications including safety assessment of radioactive waste repositories. Numerical solution of poroelasticity problems discretized with finite volume -- virtual element scheme leads to systems of algebraic equations, which may be solved simultaneously or iteratively. In this work, parallel scalability of the monolithic strategy and of the fixed-strain splitting strategy is examined, which depends mostly on linear solver performance. It was expected that splitting strategy would show better scalability due to better performance of a black-box linear solver on systems with simpler structure. However, this is not always the case.
% We would like to encourage you to list your keywords within
% the abstract section using the \keywords{...} command.
\end{abstract}
\section{Poroelasticity problem}
Modeling of coupled physical processes is important in many engineering applications, such as safety assessment of radioactive waste repositories. It is acknowledged that complex thermo-hydro-mechanical-chemical (THMC) processes should be taken into account in such modeling \cite{decovalex25}. Software package GeRa (Geomigration of Radionuclides) \cite{gera-gorny}\cite{gera-site} which is developed by INM RAS and Nuclear Safety Institute RAS already has some coupled modeling capabilities \cite{gera-thermal}\cite{gera-ruscd2018} and is now moving toward hydromechanical processes. Poroelasticity is the simplest example of such processes.

Numerical solution of coupled problems is a computationally expensive task. Arising discrete systems require efficient solution strategies, and parallel computations are a necessity. In this work, two solution strategies are tested in their scalability when a black-box linear solver is used.
\subsection{Mathematical formulation}
This work is restricted to the simplest case of elastic media filled with water only. Following theory introduced by Biot \cite{biot}, the following equations are considered:
\begin{equation}\label{eq:flow}
  s_{stor}\frac{\partial h}{\partial t} + \nabla\cdot
  \mathbf{q}
  + \alpha \nabla\cdot{\frac{\partial\mathbf{u}}{\partial t}} = Q,
\end{equation}
\begin{equation}\label{eq:mech}
  \nabla\cdot\left(\mathbf{\sigma}- \alpha P \mathbf{I}\right) = \mathbf{f}.
\end{equation}
Equation \eqref{eq:flow} represents water mass conservation taking into account porous medium deformation. Equation \eqref{eq:mech} represents mechanical equilibrium in porous medium in presence of water pressure and external forces. Here $h$ is water head, $s_{stor}$ is the specific storage coefficient, $Q$ is specific sink and source term, $\mathbf{\sigma}$ is the stress tensor, $\mathbf{f}$ is the external force vector, water pressure $P$ is related to water head $h$ as $P = \rho g \left(h - z\right)$; $\alpha$ is the Biot coefficient, which is equal to 1 in this work.  

The following constitutive relationships complete the equations: Darcy law
\begin{equation}
\mathbf{q} = -\mathbf{K}\nabla h
\end{equation}
and generalized Hooke's law:
\begin{equation}
\sigma = \mathbf{C}\varepsilon = \mathbf{C}\frac{\nabla\mathbf{u} + \left(\nabla\mathbf{u}\right)^T}{2}.
\end{equation}

Here $\mathbf{q}$ is the water flux, $\mathbf{K}$ is the hydraulic conductivity tensor, a 3$\times$3 s.p.d. matrix, $\mathbf{C}$ is the stiffness tensor, $\varepsilon$ is the strain tensor and $\mathbf{u}$ is the displacement vector.

Water head $h$ and solid displacement $\mathbf{u}$ are the primary variables.

The system is closed with initial and boundary conditions. The following boundary conditions are available:
\begin{itemize}
	\item specified head $h$ or normal flux $\mathbf{q}\cdot\mathbf{n}$ for flow;
	\item specified displacement $\mathbf{u}$, traction $\sigma\cdot\mathbf{n}$ or roller boundary condition with zero normal displacement for mechanics.
\end{itemize}

\subsection{Discretization}
Subsurface flow modeling is a well-established technology in GeRa and uses the finite volume method (FVM). Choice of the discretization method for elasticity was guided by following criteria: (a) applicability on general grids, (b) ability to work with arbitrary tensor $\mathbf{C}$, (c) sufficient history of application in multiphysics. Criterion (a) makes use of traditional finite element method (FEM) problematic, since meshes for subsurface domains can contain cells which are general polyhedra. While FVM for geomechanics exists \cite{nordbot-mpsa}\cite{terekhov-elastic} and is applied in poroelastic case \cite{terekhov-poroelastic} and more complex ones \cite{nordbot-mpfa-mpsa}, it is still somewhat new and not so widely used option. Recent developments include FVM scheme achieving improved robustness by avoiding decoupling into subproblems and introducing stable approximation of vector fluxes \cite{terekhov-poroelastic-2022}.  

For discretization of elasticity equations in GeRa the virtual element method (VEM) \cite{vem} was ultimately chosen. VEM, applied to elasticity equation \cite{vem-gain}, can handle cells which are non-convex and degenerate, its simplest version uses only nodal unknowns and is similar to FEM with piece-wise linear functions. An important feature of VEM is existence of FVM-VEM scheme for poroelasticity with proved properties \cite{coulet-coupled}. A drawback of this scheme is the use of simplest FVM option, the linear two-point flux approximation (TPFA) which is inconsistent in general case. In this work, TPFA is replaced with a multi-point flux approximation, MPFA-O scheme \cite{mpfa-o}. This scheme gives reasonable solutions on a wider class of grids, capturing media anisotropy, but is not monotone which is important for more complex physical processes.

Discretization in time uses first-order backward Euler scheme, which results in a system of linear equations at each time step.

\section{Solution strategies for discrete systems}
\begin{figure}
	\centering
	\includegraphics[width=0.4\textwidth]{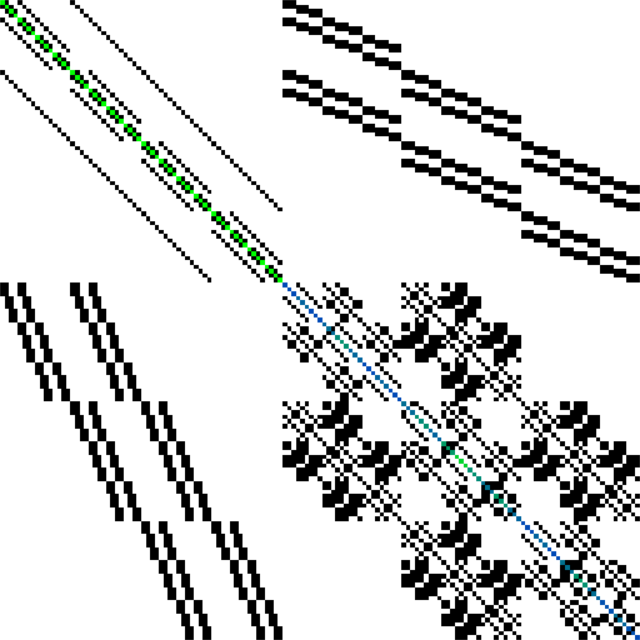}
	\caption{Matrix pattern for a 4$\times$4$\times$4 cubic grid}\label{pic:matrix}
\end{figure}

The system of discrete equations has the form
\begin{equation}\label{eq:system}
\begin{bmatrix}
A_F & A_{FM}\\
A_{MF} & A_{M}
\end{bmatrix}
\cdot
\begin{bmatrix}
h\\
\mathbf{u}
\end{bmatrix}
=
\begin{bmatrix}
b_F\\
b_M
\end{bmatrix},
\end{equation}
where subscripts $F$ and $M$ denote parts related to flow and mechanics subproblems, respectively. Here $h$ and $\mathbf{u}$ denote vectors of discrete unknowns on a given time step. The system matrix has block form with square block $A_F$ representing FVM discretization of equation \eqref{eq:flow}, square block $A_M$ representing VEM discretization of equation \eqref{eq:mech} and off-diagonal blocks $A_{FM}$ and $A_{MF}$ representing coupling terms discretized with VEM (example of matrix $A$ is depicted at figure \ref{pic:matrix}). Right-hand side terms $b_F$ and $b_M$ contain contributions from boundary conditions, source and force terms and previous time step values.

In multiphysics applications, different approaches to solution of coupled problems exist. 

\subsection{Monolithic strategy}
An intuitive approach is the monolithic strategy, in which system \eqref{eq:system} is solved as is and water head and displacement values are obtained simultaneously. This approach results in one linear system solution per time step and is unconditionally stable \cite{kim}. However, both flow and mechanics modules should be implemented in a single simulator. Moreover, complex structure of the system matrix requires robust linear solvers with some efforts being centered around specialized physics-based preconditioners \cite{erpf} and other sophisticated approaches. In case of more relevant physics like unsaturated or multiphase flow, the discrete system becomes nonlinear, and monolithic approach may lead to convergence problems of nonlinear solver, e.g. Newton method.

\subsection{Fixed-strain splitting strategy}
This solution method belongs to the class of iterative splitting schemes. Such strategies split flow and mechanics subproblems and solve each separately, iterating between then two until convergence. Notice that this splitting is not mere decomposition of matrix $A$ into blocks to solve the linear systems. The splitting involves two different solution procedures for both subproblems. With splitting strategy it possible to use tailored solvers and even separate dedicated simulators for each subproblem. One can also expect satisfactory performance from black-box linear solvers since subproblem matrices $A_F$ and $A_M$ have simpler structure compared to the full matrix $A$. In case of nonlinear flow equations, nonlinearity stays in the flow subsystem, while mechanics part remains linear and still needs to compute preconditioner only once.

Different splitting approaches are distinguished by constraints and by which subproblem is solved first. In this work, the fixed-strain approach is examined. It is the simplest splitting method where the flow subproblem is solved first. At each time step a splitting loop is executed. At each iteration of the loop the flow subproblem is solved first with displacement values fixed, then the mechanics is solved with obtained water head values. 

At each splitting iteration, residuals $r_F$, $r_M$ of flow and mechanics equations are evaluated. The splitting loop is stopped when two conditions are satisfied:
\begin{equation*}
||r_F||_2 < \varepsilon_{spl,abs} \text{~or~} ||r_F||_2 < \varepsilon_{spl,rel} \cdot ||r_F^0||_2
\end{equation*}
and
\begin{equation*}
||r_M||_2 < \varepsilon_{spl,abs} \text{~or~} ||r_M||_2 < \varepsilon_{spl,rel} \cdot ||r_M^0||_2,
\end{equation*}
where $r_*^0$ is the residual at first splitting iteration.

Fixed-strain splitting is only conditionally stable \cite{kim} and is presented here only as an example of a splitting solution strategy.

\section{Numerical experiments}
\subsection{Implementation details}
GeRa is based on INMOST \cite{inmost-site}\cite{inmost-book}, a platform for parallel computing. INMOST provides tools for mesh handling, assembly of systems via automatic differentiation as well as variety of linear solvers. Multiphysics tools of INMOST are able to switch submodels on and off in the global coupled model, which allows for easy implementation of different splitting strategies with minimal code modification.

The idea of this work is to use a black-box solver with minimal parameter tuning and compare performance of full and splitting strategies in this case.
INMOST internal solver \texttt{Inner\_MPTILUC} is used. It is a robust solver which has been successfully used in GeRa, including parallel computations in nonlinear problems with highly heterogeneous and anisotropic domains \cite{cont}. The solver is based on Bi-CGSTAB with preconditioner performing second order Crout-ILU with inverse-based condition estimation and maximum product transversal reordering \cite{inmost-book}. Convergence of the solver is governed by relative and absolute tolerances which are set to $10^{-9}$ and $10^{-12}$, respectively. Other parameters are set to default values except drop tolerance.

\subsection{Problem A: faulted reservoir}
\begin{figure}
	\centering
	\includegraphics[width=0.6\textwidth]{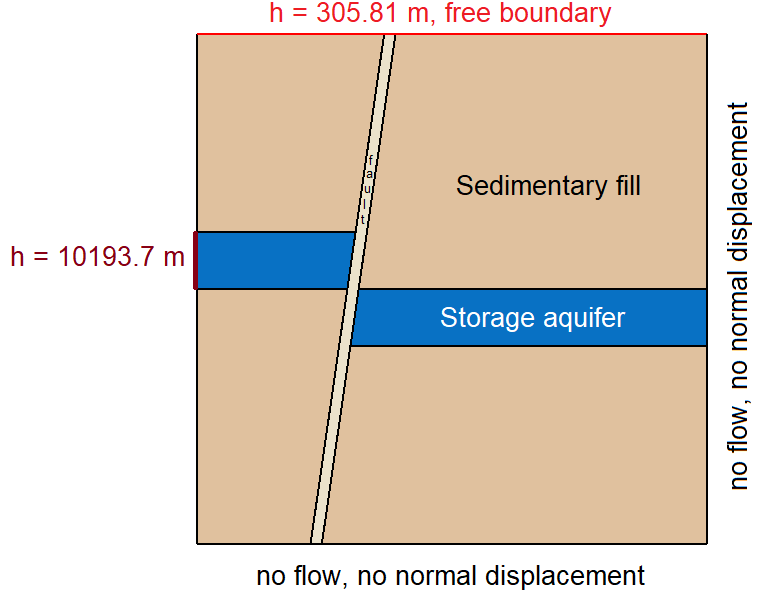}
	\caption{Problem A setup. XZ cross-section}\label{pic:a_model_setup}
\end{figure}

Problem A is a model problem similar to one presented in \cite{ogs}. It describes coupled water flow and elastic deformation in a faulted reservoir. The domain is cube 900 m $\times$ 900 m $\times$ 900 m composed of three layers. One 100 m thick layer located in the center is storage aquifer, other two layers are low-permeable sedimentary fill. An almost vertical fault crosses all three layers (see figure \ref{pic:a_model_setup}). The fault is modeled as another porous medium with increased hydraulic conductivity. All three media are isotropic, their hydraulic conductivity is characterized by a single value $K$, while stiffness tensor is completely defined by Young's modulus $E$ and Poisson ratio $\nu$. Media parameters are listed in table \ref{tab:a_media_params}.

At the top boundary, a constant water head value $h = 305.81$ m is set. %, which corresponds to pressure value $P = 3$ MPa. 
At the left side of storage aquifer, a constant water head value $h = 10193.7$ is set. %, which corresponds to pressure value $P = 100$ MPa. 
All other boundaries have zero normal flux conditions. The top boundary is free (zero traction BC), all other boundaries are sliding planes with fixed zero normal displacement. Initial water head is set constant at $h = 305.81$ m.

Simulation time covers $4\cdot 10^9$ s ($\approx 127$ years) with 4 time steps for both monolithic and fixed-strain strategies. For the linear solver, drop tolerance of 0.1 was set.

\begin{table}
	\caption{Media parameters for the problem A}\label{tab:a_media_params}
	\begin{center}
		\begin{tabular}{l | r| r| r| r}
			\hline
			Media & $K$, m/s & $s_{stor}$, 1/m & $E$, MPa & $\nu$\\
			\hline
	Storage aquifer  & 1.5$\cdot 10^{-10}$ & 8.20116$\cdot 10^{-7}$ & 14400 & 0.2\\
	Sedimentary fill &   2$\cdot 10^{-13}$ & 8.46603$\cdot 10^{-7}$ & 29400 & 0.12\\
	Fault            & 1.5$\cdot 10^{-9}$  & 1.92276$\cdot 10^{-6}$ & 14400 & 0.2
		\end{tabular}
	\end{center}
\end{table}

%\subsection{Results}
A series of tests was conducted on triangular prismatic mesh of 672300 cells and 356356 nodes, which makes total number of unknowns 1741368. Computations were carried out on INM RAS cluster \cite{cluster} on up to 100 computational cores.

\begin{figure}
	\centering
	\includegraphics[width=0.4\textwidth]{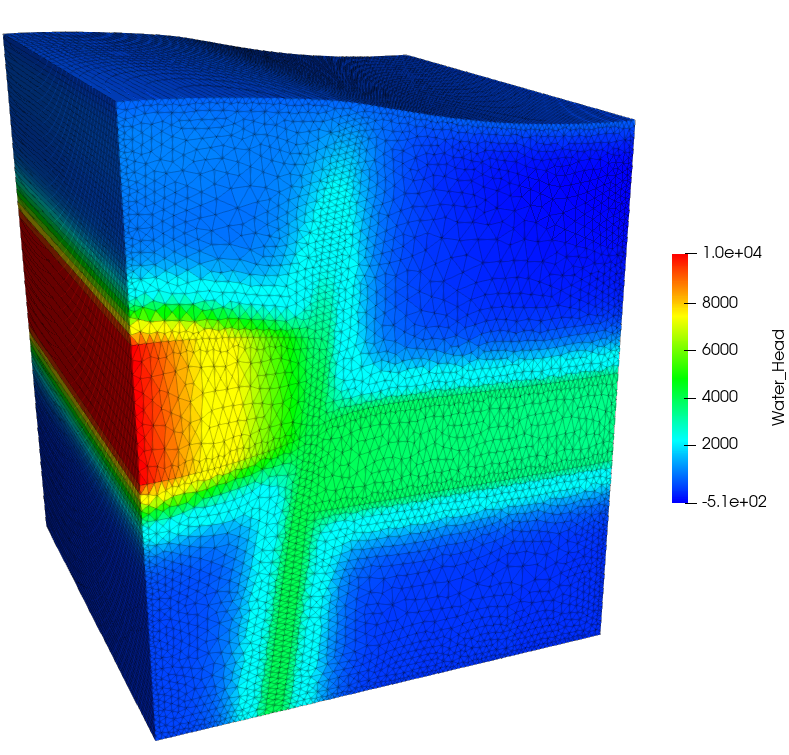}
	\includegraphics[width=0.4\textwidth]{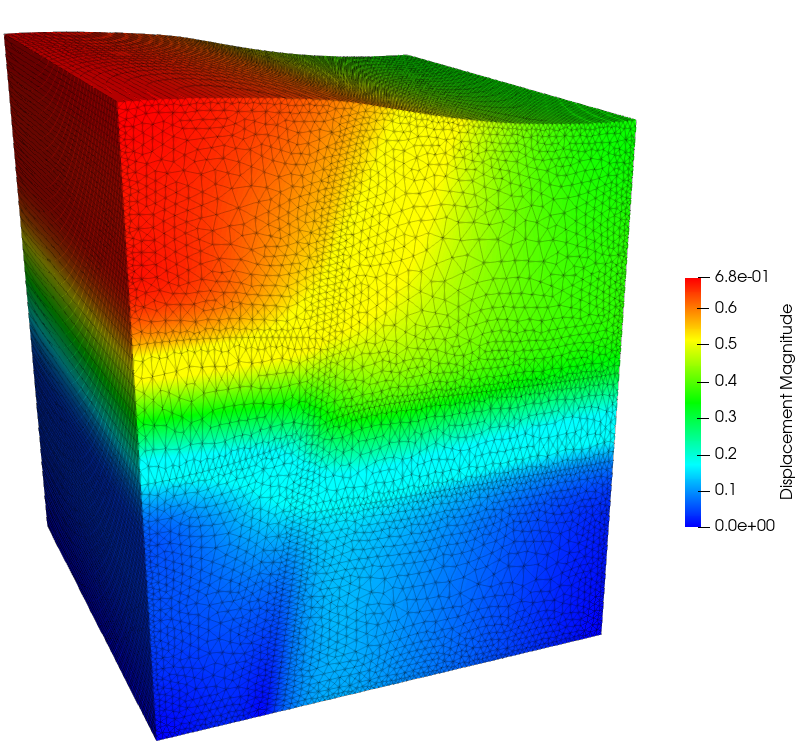}
	\includegraphics[width=0.4\textwidth]{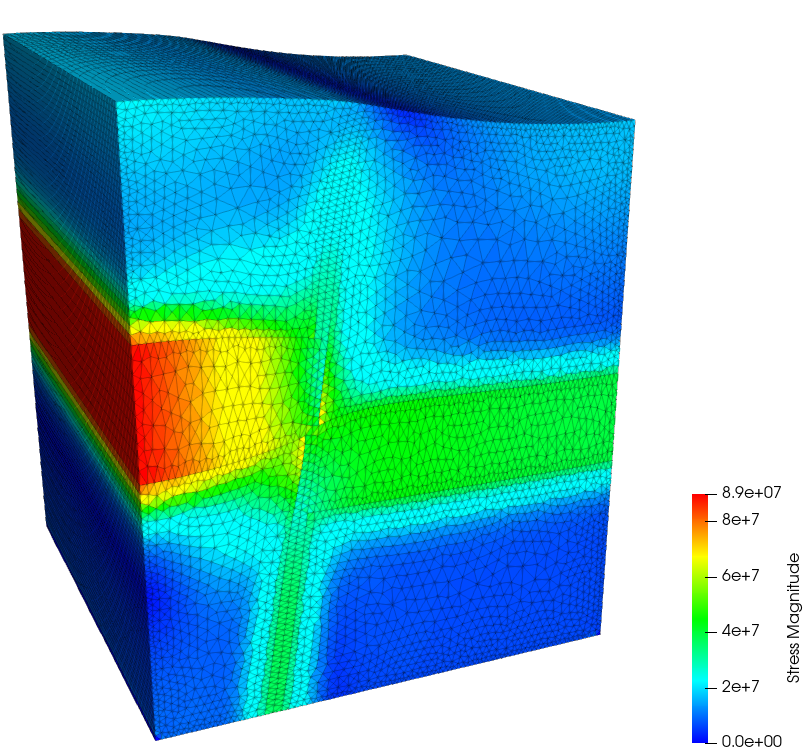}
	\caption{Problem A: water head, displacement and stress tensor in the domain, mesh distorted by displacement magnified 300 times}\label{pic:model_results}
\end{figure}

Computed water head and displacement and stress tensor magnitudes are depicted in figure \ref{pic:model_results}. Water head builds up in high-permeable aquifer and fault, resulting in stress changes which uplift surface of the domain.

\begin{figure}
	\centering
	\includegraphics[width=0.6\textwidth]{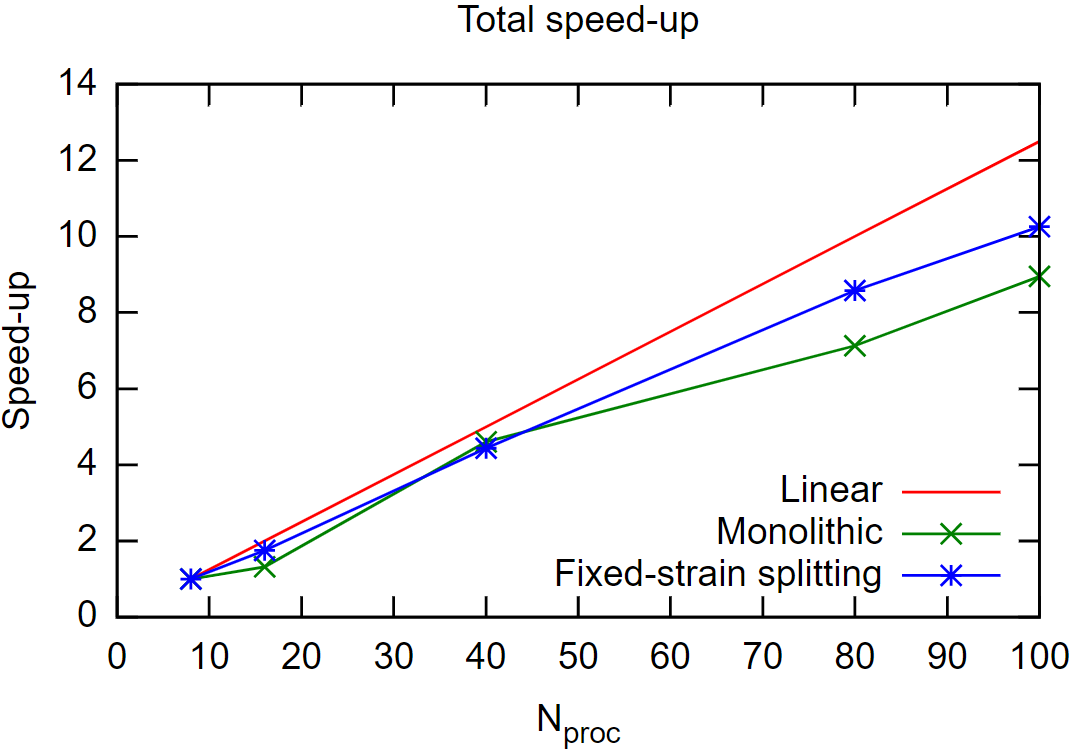}
	\caption{Problem A: total speed-up}\label{pic:speedup}
\end{figure}

\begin{figure}
	\centering
	\includegraphics[width=0.4\textwidth]{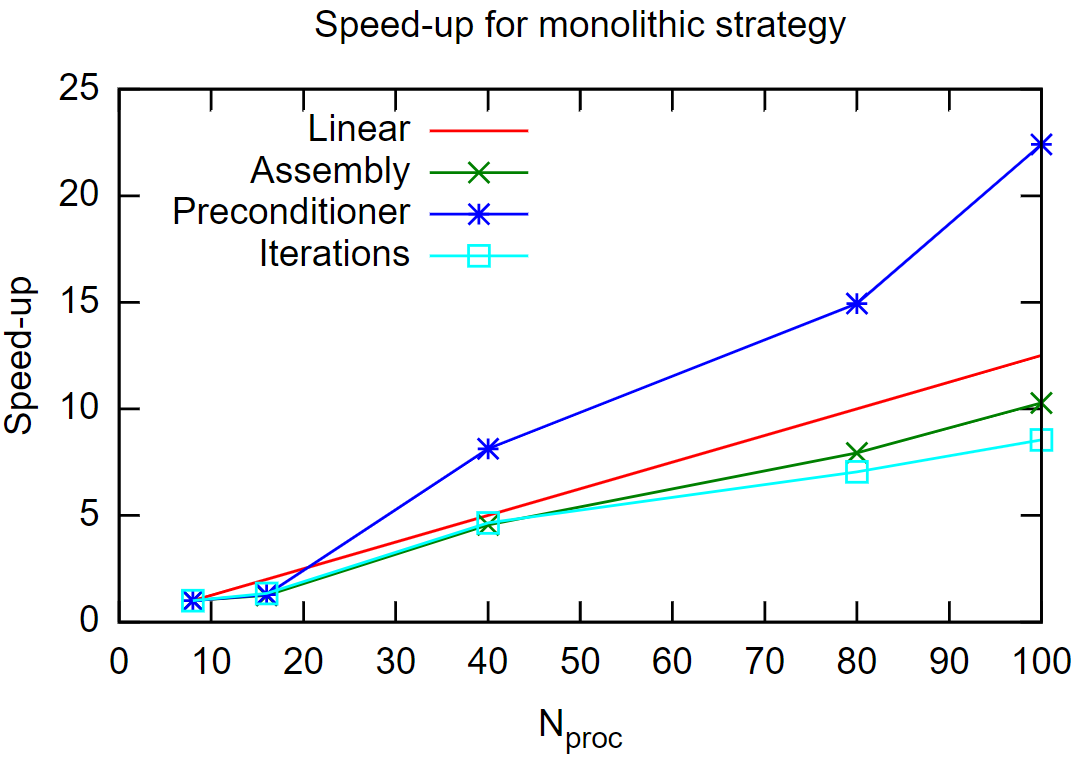}
	\includegraphics[width=0.4\textwidth]{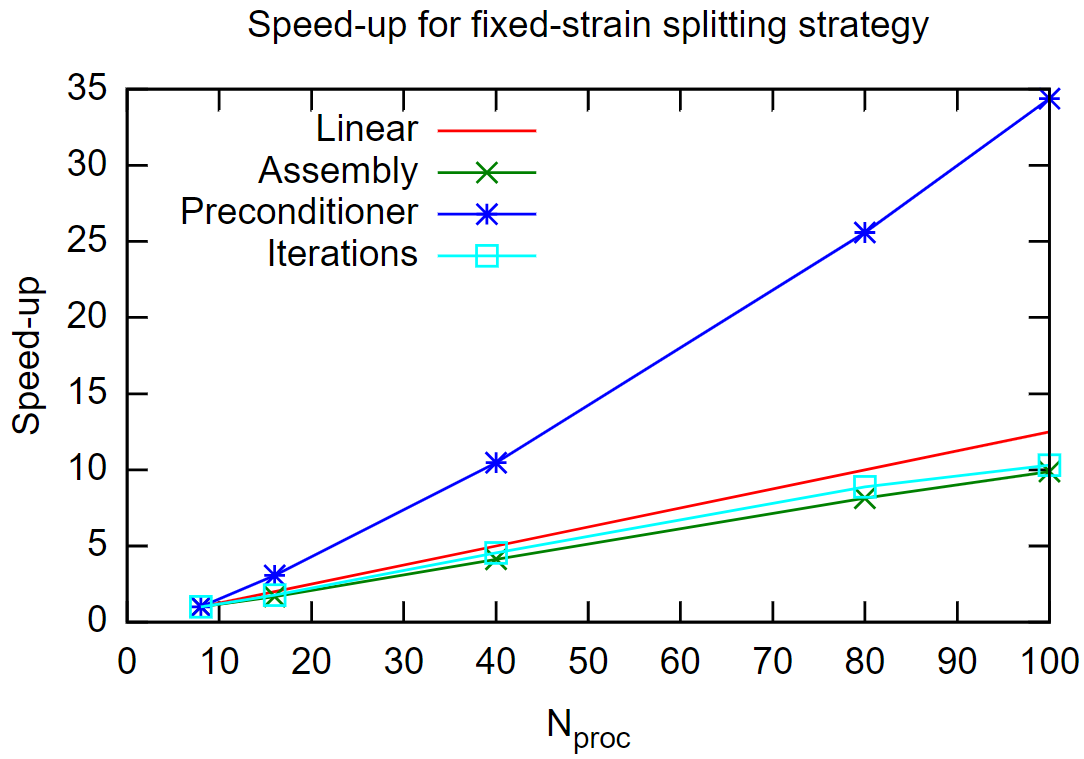}
	\caption{Problem A: detailed speed-up}\label{pic:speedup_1}
\end{figure}

\begin{figure}
	\centering
	\includegraphics[width=0.6\textwidth]{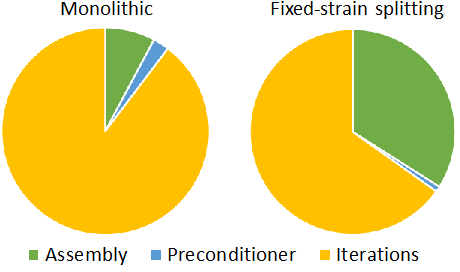}
	\caption{Problem A: time distributions for different computational stages}\label{pic:a_profiling}
\end{figure}

Measured computation times for monolithic and fixed-strain splitting methods are presented in tables \ref{tab:cluster_full} and \ref{tab:cluster_split}. Total speed-up is presented at figure \ref{pic:speedup}, while speed-ups for different computational stages for both strategies are presented in figure \ref{pic:speedup_1}. Profiling results are presented in figure \ref{pic:a_profiling}. Results show that fixed-strain strategy takes more time than monolithic one due to relatively large numbers (13 to 14) of splitting iterations at time steps. The fixed-strain strategy, however, scales better due to better scaling of iterations time, which is the most time-consuming part. Sublinear scaling of assembly procedure is explained by non-ideal partitioning by INMOST internal partitioner \texttt{Inner\_RCM}, which is based on reverse Cuthill--McKee algorithm.

\begin{table}
	\caption{Problem A: results of cluster computations, monolithic strategy}\label{tab:cluster_full}
	\begin{center}
		\begin{tabular}{l | r| r| r| r| r}
	\hline
$N_{proc}$ & $T_{total}$, s & $T_{assmbl}$, s & $T_{precond}$, s & $T_{iter}$, s & \#lin.it\\ \hline
8          &  2103          &   184           &  130             &  1761         &   1673   \\
16         &  1592          &   148           &  101             &  1318         &   1637   \\
40         &   457          &    40.3         &   16             &   379         &   2298   \\
80         &   295          &    23.2         &    8.7           &   250         &   2522   \\
100        &   235          &    17.9         &    5.8           &   206         &   2639   \\
	
\end{tabular}
	\end{center}
\end{table}

\begin{table}
	\caption{Problem A: results of cluster computations, fixed-strain splitting strategy}\label{tab:cluster_split}
	\begin{center}
		\begin{tabular}{l | r| r| r| r| r}
			\hline
$N_{proc}$ & $T_{total}$, s & $T_{assmbl}$, s & $T_{precond}$, s & $T_{iter}$, s & \#lin.it\\ \hline
8          &  7684          &  2479           &  219.5           &  4961         &  15162   \\
16         &  4390          &  1483           &   71.5           &  2817         &  18121   \\
40         &  1733          &   601.2         &   21.1           &  1091         &  19575   \\
80         &   896          &   305.6         &   8.6            &   558.8       &  21099   \\
100        &   749          &   251.7         &   6.4            &   480.8       &  21395   \\
	 
		\end{tabular}
	\end{center}
\end{table}

\subsection{Problem B: real-life domain with synthetic elastic parameters}
In this problem, a part of a real site is considered. The part is a quadrilateral cut in XY-plane from a domain with 9 geological layers and 11 different media. The media have are anisotropic with their hydraulic conductivity being a diagonal tensor with values ranging from $1.2\cdot10^{-12}$ to $2\cdot10^{-5}$ m/s. Specific storage coefficient varies from $10^{-6}$ to $10^{-5}$ m$^{-1}$. Elastic parameters are not known and are set constant for all media: $E = $ 10000 MPa, $\nu = 0.2$. One corner of the model has prescribed water head value $h = 1000$ m at 8th layer, which imitates injection in that layer. On the opposite side of the domain, constant water head value of $h = 100$ m is set. Other boundaries are impermeable. Bottom boundary is fixed ($\mathbf{u} = 0$), top boundary is free ($\sigma\cdot\mathbf{n} = 0$), side boundaries have roller boundary conditions (zero normal displacement $\mathbf{u}\cdot\mathbf{n} = 0$). Simulation starts with constant initial water head value $h = 100$ m and covers $2\cdot10^8$ s ($\approx$ 6.3 years).

\begin{figure}
	\centering
	\includegraphics[width=0.45\textwidth]{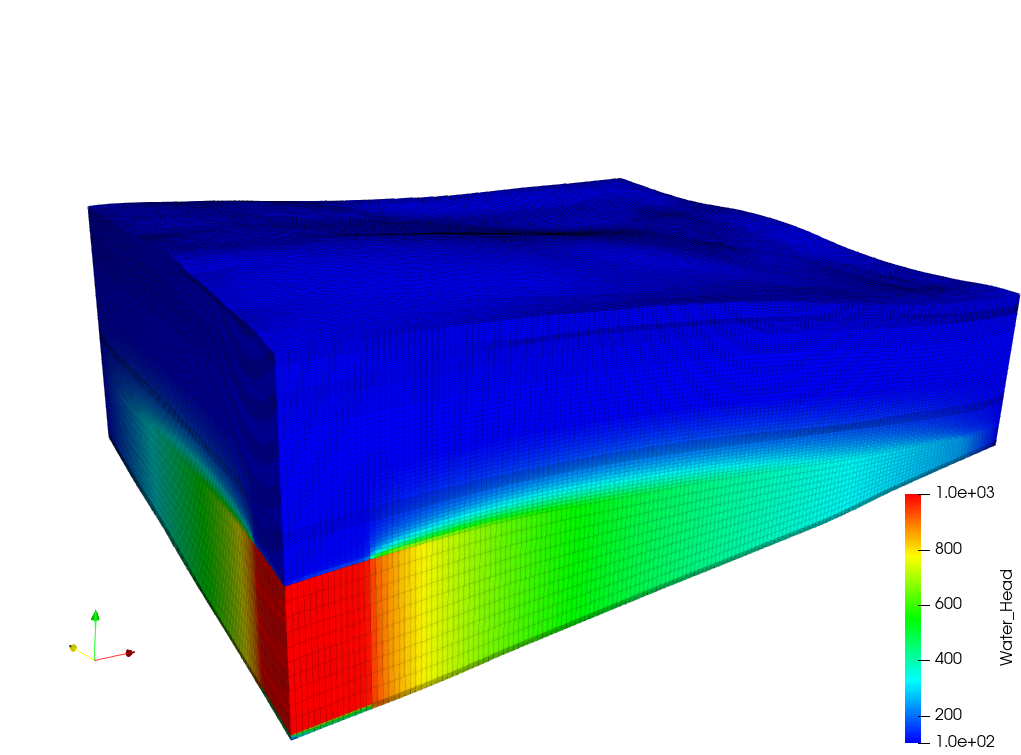}
	\includegraphics[width=0.45\textwidth]{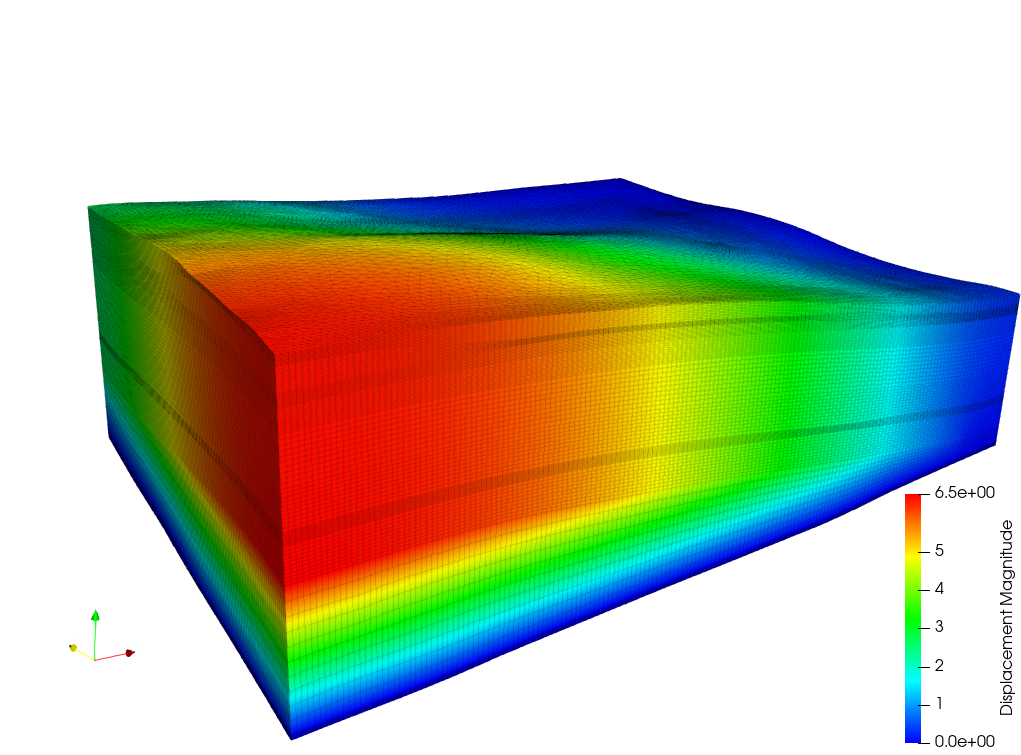}
	\caption{Problem B: water head and displacement in the domain, mesh distorted by displacement magnified 30 times}\label{pic:b_results}
\end{figure}

\begin{figure}
	\centering
	\includegraphics[width=0.6\textwidth]{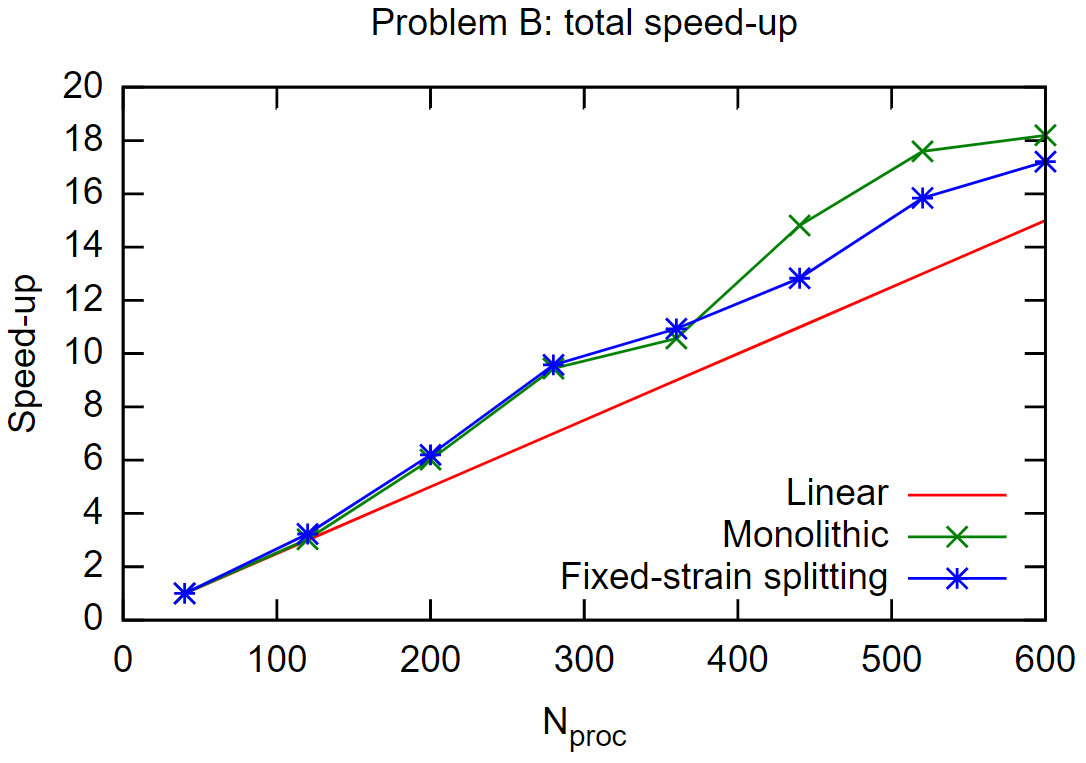}
	\caption{Problem B: total speed-up}\label{pic:b_speedup}
\end{figure}

\begin{figure}
	\centering
	\includegraphics[width=0.4\textwidth]{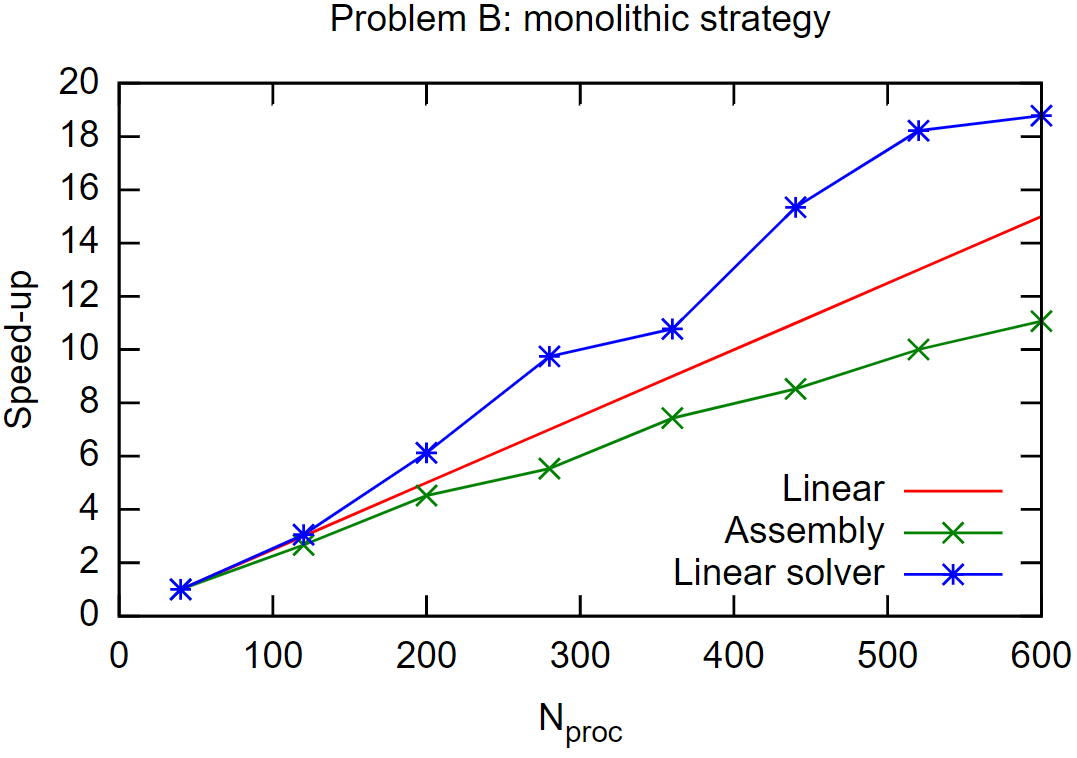}
	\includegraphics[width=0.4\textwidth]{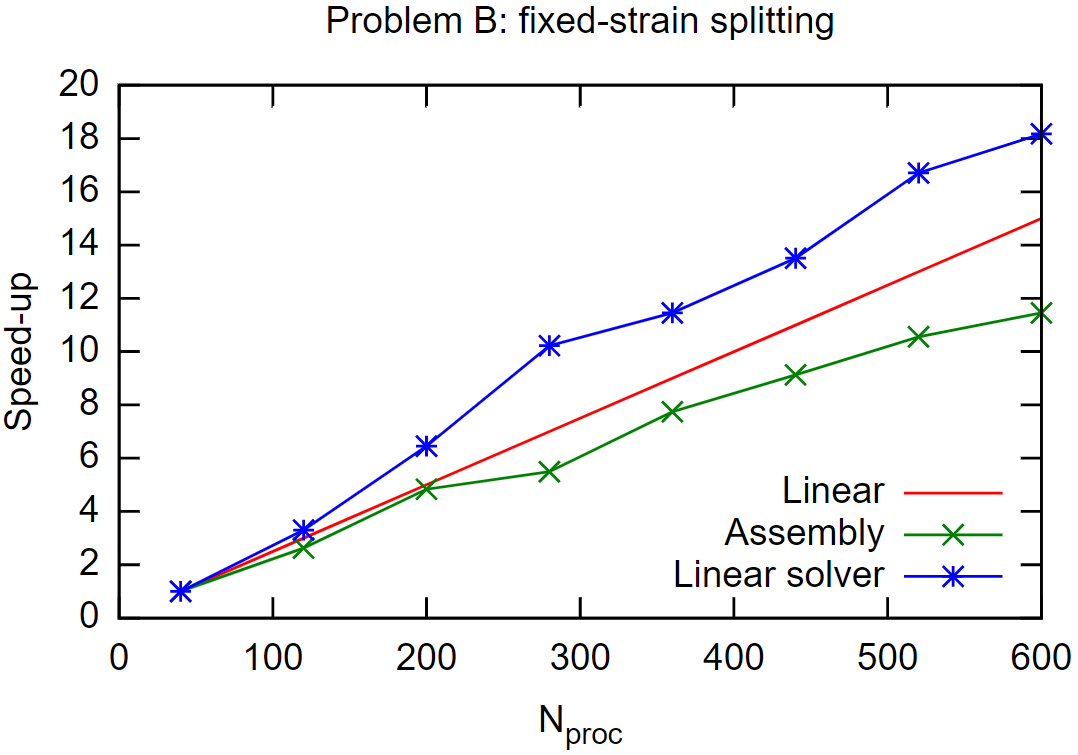}
	\caption{Problem B: speed-up for assembly and linear solver}\label{pic:b_speedup_stages}
\end{figure}

\begin{figure}
	\centering
	\includegraphics[width=0.6\textwidth]{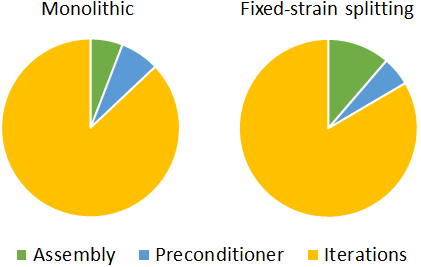}
	\caption{Problem B: time distributions for different computational stages}\label{pic:b_profiling}
\end{figure}

Triangular prismatic mesh of 2205400 cells and 1142714 (total number of unknowns is 5461942) was constructed in the domain. Computations were performed on INM cluster using 40 -- 600 cores. In order to obtain better balanced partitioning of the mesh, the ParMETIS \cite{parmetis} partitioner was used, interface to which is provided by INMOST. Calculated water head and displacement distributions are depicted in figure \ref{pic:b_results}. Time measurements and linear iterations count are presented in tables \ref{tab:b_full} and \ref{tab:b_split}. Total speed-up is presented in figure \ref{pic:b_speedup} and shows superlinear scaling with monolithic strategy reaching slightly larger speed-up. In order to better understand this scaling behavior of the two strategies, speed-up for both assembly and linear solver are presented in figure \ref{pic:b_speedup_stages} and profiling results are presented in figure \ref{pic:b_profiling}. It can be seen that superlinear scaling is caused by linear solver, namely, scaling of the preconditioner. Such scaling is caused by the fact that drop tolerance, which is left at default value of $10^{-5}$, which makes MPTILUC preconditioner close enough to performing full LU decomposition of matrix diagonal sub-blocks. At the same time, assembly scales sublinearly, which is explained again by non-ideal partitioning of the mesh. Since profiling shows that assembly takes large percentage of time in the fixed-strain splitting, this strategy reaches lower overall speed-up.

\begin{table}
	\caption{Problem B: results of cluster computations, monolithic strategy}\label{tab:b_full}
	\begin{center}
		\begin{tabular}{l | r| r| r| r| r}
			\hline
$N_{proc}$ & $T_{total}$, s & $T_{assmbl}$, s & $T_{precond}$, s & $T_{iter}$, s & \#lin.it\\ \hline
40         &  1490          &    52.0         &  688             &   731.5       &    470   \\
120        &   490.5        &    19.5         &   58.0           &   406.7       &    680   \\
200        &   247.5        &    11.5         &   51.1           &   180.5       &    742   \\
280        &   157.8        &     9.4         &   16.2           &   129.5       &    779   \\
360        &   141.0        &     7.0         &   24.6           &   107.1       &    902   \\
440        &   100.6        &     6.1         &   10.1           &    82.4       &    821   \\
520        &    84.7        &     5.2         &    6.3           &    71.6       &    840   \\
600        &    81.9        &     4.7         &    5.7           &   69.9        &    993   \\

		\end{tabular}
	\end{center}
\end{table}

\begin{table}
	\caption{Problem B: results of cluster computations, fixed-strain splitting strategy}\label{tab:b_split}
	\begin{center}
		\begin{tabular}{l | r| r| r| r| r}
			\hline
$N_{proc}$ & $T_{total}$, s & $T_{assmbl}$, s & $T_{precond}$, s & $T_{iter}$, s & \#lin.it\\ \hline
 40        &  1101          &    81.3         &  475.3           &  537.4        &    728   \\
120        &   340.3        &    31.0         &  28.9            &  278.1        &   1079   \\
200        &   177.4        &    16.8         &  32.0            &  125.0        &   1101   \\
280        &   114.9        &    14.8         &    8.6           &   90.5        &   1175   \\
360        &   100.7        &    10.5         &   15.2           &   73.2        &   1290   \\
440        &    85.8        &     8.9         &    5.5           &   69.5        &   1408   \\
520        &    69.5        &     7.7         &    3.5           &   57.1        &   1399   \\
600        &    64.0        &     7.1         &    3.3           &   52.4        &   1503   \\			
		\end{tabular}
	\end{center}
\end{table}

\section{Conclusion}
Two solution strategies for finite volume -- virtual element discretizations of coupled poroelastic problems were tested in their parallel performance. Implemented using INMOST numerical platform, the model used a robust general-purpose linear solver. Calculations for two problems on meshes with millions of unknowns and up to 600 computational cores gave mixed results. Fixed-strain splitting strategy was expected to scale better due to performance of the linear solver, as separate systems for flow and mechanics have simpler structure. This, however, was not always observed. This is explained by the following reasons. First, performance of the linear solver is mostly determined by preconditioner, which exhibits superlinear scalability depending on the problem. Second, non-ideal performance of mesh partitioner results in sublinear scalability of systems assembly process, which affects splitting strategy more as assembly takes larger fraction of time in this case. Overall, there was no clear answer as to which strategy scales better. Partly this may imply that the linear solver used can successfully handle fully coupled multiphysical systems. However, with additional partitioner tuning to better balance mesh between processes, scalability of assembly may approach linear rate, resulting in better overall scalability of splitting strategy. 


\end{document}